\documentclass[11pt,epsf]{article}
\usepackage{graphicx}
\usepackage[italian, english]{babel}	% per scrivere in italiano... utile per gli accenti
\usepackage{amsmath}
\usepackage[most]{tcolorbox}
\usepackage{anyfontsize}
\usepackage{courier} %% Sets font for listing as Courier.
\usepackage{listings, xcolor}
\newtheorem{theorem}{Theorem}[section]
\newtheorem{lemma}[theorem]{Lemma}
\newtheorem{definition}[theorem]{Definition}

\newenvironment{thm}{\begin{theorem}}{\end{theorem}}

\newenvironment{proof}[1]{{\noindent\bf Proof:~}#1}{\qed}

\def\whitebox{{\hbox{\hskip 1pt
         \vrule height 6pt depth 1.5pt
         \lower 1.5pt\vbox to 7.5pt{\hrule width
                   3.2pt\vfill\hrule width 3.2pt}%
         \vrule height 6pt depth 1.5pt
         \hskip 1pt } }}

\def\qed{\ifhmode\allowbreak\else\nobreak\fi\hfill\quad\nobreak
     \whitebox\medbreak}

\definecolor{darkgreen}{RGB}{0,127,0}

\newcommand{\cfor}{{\text{ for }}}
\definecolor{block-gray}{gray}{0.85}

\begin{document}
\title{A simple proof for \\generalized Fibonacci numbers with dying rabbits\footnote{
A full version of this paper appeared as {\em A combinatorial proof for the Fibonacci dying rabbits problem}, in the Bulletin of the Institute of Combinatorics and
its Applications, Volume 103 (2025), pp. 25--36~\cite{D2025}. We refer the reader to the full version for numerical examples.
In this version we have added reference~\cite{F2011} and some remarks about it (see the ``Previous work'' section).}}
\author{Roberto De Prisco\\
{\small Dipartimento di Informatica}\\
{\small University of Salerno, Italy}\\
{\small\tt robdep@unisa.it}}
\date{\today}
\maketitle

\begin{abstract}
We consider the generalized Fibonacci counting problem with rabbits that become
fertile at age $f$ and die at age $d$, with $1\le f\le d$, and $d$ finite or infinite. 
We provide a simple proof, based exclusively on a counting argumentation, for a recurrence relation.
Denoting with $F_n$ the number of rabbits at generation $n$, with the initial condition $F_1=1$, we have that

\begin{equation*}
F_n = 
\begin{cases}
1,                                         & \cfor 2\le n\le f,\\
F_{n-1}+F_{n-f},                   & \cfor f+1\le n\le d,\\
F_{n-1}+F_{n-f}-1,              & \cfor n=d+1,\\
F_{n-1}+F_{n-f}-F_{n-d-1}, & \cfor n\ge d+2.\\
\end{cases}
\end{equation*}
\noindent
This formula reduces to the classical Fibonacci case when $f=2$ and $d=\infty$.
\end{abstract}

\section{Introduction}

Leonardo Bonacci, better known as Fibonacci, the $13^{th}$ century Italian mathematician, among
his many contributions to the field, considered
the well known rabbits counting problem, that resulted in the so-called Fibonacci sequence
of integers. The problem is the following. A population of pair of rabbits, starting with one newborn pair, grows
in each {\em generation} with every fertile pair of rabbits giving birth to a new pair of rabbits.
A pair of rabbits becomes fertile at age 2, that is, it does not proliferate in the generation in which it
is born but starts proliferate in the next generation, and the newborns are added the subsequent
generation. The problem is that of counting the number $F_n$ of pairs of rabbits at every generation.
The initial condition gives $F_1=1$. The unique pair of rabbits is not fertile in the first generation and thus $F_2=1$.
In the second generation the pair proliferates giving birth to a newborn pair of rabbits and thus $F_3=2$.
For the next generation only the initial pair of rabbits proliferates and thus $F_4=3$. For the subsequent one
there are two pairs of rabbits that proliferates and thus $F_5=5$. Proceeding in this way one gets the
so-called Fibonacci sequence
$$1,1,2,3,5,8,13,21,34,55,89,\ldots.$$
The well known recurrence relation
\begin{equation}\label{eq:fibonacci}
F_n=F_{n-1}+F_{n-2},
\end{equation}
gives an easy way to compute the Fibonacci sequence, for any $n\ge 3$, starting from the initial conditions
$F_1=1$ and $F_2=1$. One can extend the sequence by adding $F_0=0$, and the formula
is still valid for $n\ge 2$, with the initial conditions $F_0=0$ and $F_1=1$.

It is well know that the Fibonacci sequence has many interesting properties that have lead to its fame.
We refer to the many papers and books available in the literature for the reader interested in these properties.

\medskip

A number of generalizations have been considered. For example the Lucas sequence, studied
in the 19$^{th}$ century by the French mathematician François Édouard Anatole Lucas, is
obtained by changing the initial conditions to $F_0=2$ and $F_1=1$.
The Gibonacci sequence is a further generalization for which $F_0$ and $F_1$ can be arbitrary;
thus Fibonacci and Lucas numbers are special cases of Gibonacci numbers. 
The Padovan sequence is defined as $F_n=F_{n-2}+F_{n-3}$, with initial numbers $F_0=F_1=F_2=1$ and
has properties similar to the Fibonacci sequence (and it actually corresponds to  one case of the generalization that we will consider).
The Tribonacci numbers are obtained adding the 3 previous numbers, while the Tetranacci numbers
are obtained adding the previous 4 numbers in the sequence. More in general sequences
obtained by adding the previous $k$ elements, that is, with $F_n=\sum_{i=1}^k F_{n-i}$, have been considered.
Another generalization is called $k$-Fibonacci and the numbers are defined as $F_n=k F_{n-1}+F_{n-2}$.

There are plenty of research papers about these generalizations and several books have been written about
the Fibonacci numbers; 
we refer the reader to~\cite{book:K2001,book:V1989} for more information.

\medskip

In this paper we consider the generalization in which, in the original problem, cast as a counting 
problem of a growing population of rabbits, the rabbits become fertile after some number $f$ of
generations and at some point, after $d$ generations, they die. This problem is also called the
dying rabbits problem and has been studied in several papers~\cite{A1963-1,A1963-4,C1964,H1968,HL1969,O2009}. 

\medskip
\noindent{\bf Contribution of this paper.} In this paper we provide a recurrence relation that
gives the $n^{th}$ generalized Fibonacci number as a function of 2 or 3 previous numbers. 
With the initial condition $F_1=1$, we have that

\begin{equation*}
F_n = 
\begin{cases}
1,                                         & \cfor 2\le n\le f \makebox{\ \ \ \ \ (Case 1),}\\
F_{n-1}+F_{n-f},                   & \cfor f+1\le n\le d \makebox{\ \ \ \ \ (Case 2),}\\
F_{n-1}+F_{n-f}-1,              & \cfor n=d+1  \makebox{\ \ \ \ \ (Case 3), and}\\
F_{n-1}+F_{n-f}-F_{n-d-1}, & \cfor n\ge d+2  \makebox{\ \ \ \ \ (Case 4)}.\\
\end{cases}
\end{equation*}

\medskip

\noindent
This formula clearly generalizes Equation~(\ref{eq:fibonacci}). Indeed the Fibonacci sequence
is the case for which $f=2$  and $d=\infty$ and for this choice of the two parameters, we only have 
cases 1 and 2. For $n=2$, we have case 1, which gives $F_2=1$ and for $n\ge 3$, case 2 becomes
$F_n= F_{n-1}+F_{n-2}.$

\noindent
The proof that we provide is quite simple and involves only counting arguments.

%%%%%%%%%%%
%%%%%%%%%%%
%%%%%%%%%%%
\section{Previous work}

The earliest work that we are aware of, about the Fibonacci counting problem with dying rabbits, 
are by Brother U. Alfred who, in~\cite{A1963-1}, posed the question of counting the rabbits
for the specific case of $f=2$ and $d=12$, probably thinking that it was a relatively easy counting problem
and later, in~\cite{A1963-4}, concluded that the problem did not seem that easy. 
A recurrence relation that correctly solves the case $f=2$, $d=12$ has been provided by Cohn~\cite{C1964}.
The recurrence relation matches the general formula that we provide in this paper.

In~\cite{H1968} and~\cite{HL1969} a generalization on
the breeding pattern is considered. Each pair of rabbits breeds $B_1$ new pairs in its first generation,
$B_2$ in the second and so on, with $B_0=0$, and with rabbits dying after a fixed number
of generations. The proposed solution is given in terms of generating functions, and no
explicit formula is derived.

Oller-Marcén~\cite{O2009} provides the following recurrence relation.
The notation
used in~\cite{O2009} uses $h$ for the age in which rabbits become fertile and $k$ for the generations
that rabbits live {\em after} the fertile age. Thus the correspondence with our notation is $f=h$ and $d=k+h-1$.
The recurrence relation is

\begin{equation}\label{eq:oller}
F_n = 
\begin{cases}
1,                                                       & \cfor 1\le n\le f\\
F_{n-1}+F_{n-f},                                & \cfor f < n \le d\\
\underbrace{F_{n-f}+F_{n-f-1}+\ldots+F_{n-d}}_{\makebox{{\small $d-f+1$ terms}}}, & \cfor n > d.
\end{cases}
\end{equation}
The proof of the above equation provided in~\cite{O2009}, in Proposition 9, is quite short and although the first two cases
are easy and clear, the third case, the one for $n\ge d$, is not as much clear.
Literally --- we have only adjusted
the terminology and the notation to make them conform to the ones we are using in this paper --- the proof states:
{\em ``the number of rabbits at the $n^{th}$ generation can be computed as the sum of all the preceding
rabbits except those which are not mature yet ($F_{n-j}$, with $1\le j \le f-1$) and those which 
have died ($F_{n-j}$, with $j>n-d$)''}. Although this sentence exactly identifies the right hand side of
the third case of Equation~(\ref{eq:oller}), it is not clear why it should give the total number of rabbits for
the $n^{th}$ generation. Moreover it is not true that the number of rabbits that are not mature in generation $n$ is equal
to $F_{n-1}+F_{n-2}+\ldots+F_{n-f+1}$, as this sum is the sum of all the alive rabbits in a number of generations, and thus 
non-mature rabbits would be counted multiple times. Moreover some fertile rabbits, alive in the previous generations, might have died before generation $n$. 
A similar remark, about the multiple counting, applies to the number of rabbits that have died.

The approach used in the proof that 
we provide in this paper, leading to a more compact formula, can probably be used to explain also~(\ref{eq:oller}),
which, in fact, is equivalent to the one we propose. 
In~\cite{O2009} the solution is given also
as a function of the complex roots of the characteristic polynomial.

After the publication of~\cite{D2025} we found~\cite{F2011} that actually provides the same formula that we provide in this paper.
However the proof provided in~\cite{F2011} is wrong. Literally, adjusting the notation used in~\cite{F2011} which is the same
as the one of~\cite{O2009}, for which we have $f=h$ and $d=k+h-1$, the proof (of Theorem 2.1) says:
``{\em If $n\ge d+2$, the number of rabbits at the $n^{th}$ generation can be computed as the sum of
the pairs at the preceding generation $F_{n-1}$, and those bred by the pairs which are mature at this
point, i.e. $F_{n-f}$. (It is the number of pairs of rabbits which just came into the world.) Finally,
the number of pairs of rabbits which died at this point, i.e. $F_{n-d-1}$, must be subtracted.}''

It is not true that the number of rabbits that are born for (that is, added to) generation $F_n$ is $F_{n-f}$ and
it is not true that the number of rabbits that die for (that is, will be not part of) generation $F_n$ is $F_{n-d-1}$.
If we call $newborns_n$ and $deaths_n$, as we will do in the subsequent sections of  the paper, what the above proof is saying is that
$newborns_n=F_{n-f}$ and $deaths_n=F_{n-d-1}$. But this is not true because there can be rabbits of generation $n-f$ that have died and thus
$F_{n-f}$ is an overestimate of $newborns_n$ and 
the same is true for generation $n-d-1$ and thus $F_{n-d-1}$ is an overestimate of $deaths_n$. As an example consider
generation number 20, for the case $f=3$ and $d=9$ (see the table in Appendix B of~\cite{D2025}).
There are $F_{20}=715$ rabbits, with the following ages: 
228 rabbits of age 1,
158 rabbits of age 2,
109 rabbits of age 3,
76 rabbits of age 4,
53 rabbits of age 5,
36 rabbits of age 6,
25 rabbits of age 7,
18 rabbits of age 8, and
12 rabbits of age 9.
Thus the number of fertile rabbits that proliferates for generation 21 is the sum of rabbits with
age greater or equal to 3, that is 109+76+53+36+25+18+12=329.
The number of rabbits that die and will not be part of generation 21 is exactly the number of
rabbits with age 9, that is 12.
Thus $newborns_{21}=329$ and $deaths_{21}=12$. But $F_{n-f}=F_{18}=343$ and $F_{n-d-1}=F_{11}=26$.

Incidentally, it happens to be true that  $newborns_n-deaths_n$ is equal to $F_{n-f}-F_{n-d-1}$, and this explains why the formula is true.
In the previous example we have that $329-12=317$ and $343-26=317$.
The equality $newborns_n-deaths_n=F_{n-f}-F_{n-d-1}$ is proved in this paper.

%%%%%%%%%%%
%%%%%%%%%%%
%%%%%%%%%%%
\section{The ``base equation''}

First, to ease the wording, we will talk about single rabbits instead of pairs of rabbits. A single rabbit
that proliferates it is not natural\footnote{Also the setting of the original problem is not very natural!}, but considering single rabbits instead of pairs
of rabbits does not change the underlying counting problem.
Let $F_n$ be the number of rabbits at the $n^{th}$ generation, with the initial condition $F_1=1$.
The rabbits become fertile at the $f^{th}$ generation, and die at the age of $d$, with $1\le f \le d$.

Rabbits with age $d$ first proliferate and then die. A newborn rabbit has age $1$. 
The ``step of the evolution''
is as follows: given a population $F_n$, every element of the population with age at least $f$
gives birth to a new element for the next population; then the age of each element is increased by $1$
and rabbits with age $>d$ die and thus will not be part of the next generation.

The Fibonacci sequence is the special case $f=2$ and $d=\infty$, in which the rabbits become fertile
after their first generation, that is at the second generation, and never die.

Consider the $n^{th}$ generation and define $newborns_n$ and $deaths_n$ as the number of,
respectively, newborn rabbits and deaths. Notice that $newborns_n$ is equal to the number
of fertile rabbits in the previous $(n-1)^{th}$ generation and
$deaths_n$ is equal to the number of rabbits that have exactly age $d$ in the
previous $(n-1)^{th}$ generation.
The following basic fact should be immediate

\bigskip

\noindent{\bf Base equation.}
The number of rabbits at the $n^{th}$ generation is equal to the number of rabbits in the previous
generation plus the number of newborns minus the number of deaths, that is
\begin{equation}\label{eq:base}
F_n=F_{n-1} + newborns_{n} - deaths_{n}.
\end{equation}

\bigskip

\noindent
We call the above equation the {\em base equation}.
The base equation gives, in a straightforward way, a solution, once we are able to estimate the
newborns and the deaths at each generation. 

\bigskip

Before we proceed,
let us make some observations.
The degenerate case $d<f$ would give rise to the sequence 
$$\underbrace{1,1,\ldots,1}_{d\makebox{ {\tiny times}}},0,0,0,0,0,0,0,0,0,0,\ldots\ldots$$
since the unique rabbit will die before proliferating. 
Thus it is not interesting and this is why we consider only the case $f\le d$.

The borderline case $f=d$ is simple to deal with. Indeed, in this case every rabbit
proliferates in the same generation in which it dies, that is $newborns_n=deaths_n$ (0 for the first $f$ generations, 1 thereafter). 
Hence the total number of rabbits never changes. Thus the sequence that we get is
$$1,1,1,1,1,1,1,1,1,1,\ldots\ldots.$$

\medskip

When $f<d$ estimating the number of newborns and deaths becomes trickier; however when $d=\infty$ we can easily estimates
the deaths: $deaths_n=0$.

The special case $f=1$ (and $d=\infty$), for which rabbits are immediately fertile leads to a doubling of the rabbits
at every generation. Indeed we would have $newborns_n=F_{n-1}$, leading to $F_n=2F_{n-1}$, which gives the sequence of powers of 2:
$$1,2,4,8,16,32,64,128,256,1024,\ldots\ldots.$$

For the original Fibonacci sequence, beside $d=\infty$, we have $f=2$, and also in this case it is easy to
estimate the newborns: among the previous population of size $F_{n-1}$ there are exactly $F_{n-2}$
fertile elements, since
rabbits become fertile at age $2$, and thus there are exactly $F_{n-2}$ newborns for the new generation,
that is, $newborns_n=F_{n-2}$. 
Thus for the case $f=2$, $d=\infty$, we have the well known Fibonacci's formula $F_n=F_{n-1}+F_{n-2}.$

More in general, for the case $d=\infty$ and any finite $f$, we have that the number of rabbits that are
fertile for generation $n$ is exactly $F_{n-f}$, that is, all the rabbits that were in the population $f$
generations before; indeed all these rabbits have age at least $f$ in generation $n$ and they have not died.
All the other rabbits are still too young. The number of
fertile rabbits gives the number of newborns, that is, $newborns_n=F_{n-f}$.
This results in the formula $F_n=F_{n-1}+F_{n-f}$, as stated in Equation~(\ref{eq:oller}).

%%%%%%%%%%%
%%%%%%%%%%%
%%%%%%%%%%%

However, when $d$ is finite, estimating the exact number of newborns and deaths seems trickier. 
Indeed it is not true anymore that the number of fertile rabbits for the $n^{th}$ generation is equal
to $F_{n-f}$ because some of those rabbits could have died meanwhile. Also the total number of deaths
seems more difficult to assess.

\section{Unraveling the base equation}

Keeping track only of the total number of rabbits in each generation is not helpful.
We found that estimating also the number of rabbits for each possible age,
helps in clearly defining the relations among the numbers that we get from counting the rabbits. 
In the following we first provide some basic definition and properties, then, exploiting such properties,
we prove the four cases of the proposed formula.

\subsection{Definitions and properties}

We start with the following definition for the number of rabbits with a specific age in a given generation.

\begin{definition}\label{def:F_n^x}
Define $F_n^x$, for $x=1,2,\ldots,d$ as the number of rabbits of age $x$ at (the beginning of) generation $n$.
\end{definition}

\smallskip

\noindent
In the following, we  study the relation between all the $F_n^x$, for any $n$ and $x$, and the total number of rabbits in each generation,
that is $F_n$, for any $n$. 

\medskip

\noindent
We start with an obvious relation which follows directly from the definition of the $F_n^x$.
\begin{equation}\label{eq:Fn=sum_of_row}
F_n = F_n^1 + F_n^2 + \ldots F_n^{d-1} + F_n^d.
\end{equation}

\begin{lemma}\label{lem:goback}
For any $x\le\min\{d,n\}$, we have that
$F_n^{x}=F_{n-1}^{x-1}=F_{n-2}^{x-2}=\ldots=F_{n-x+1}^1.$
\end{lemma}

\begin{proof}
The age of the rabbits increases by 1 at each generation. Thus, if there are $F_n^x$ rabbits (of age $x$)
at generation $n$, there must have been the same number $F_{n-1}^{x-1}$ (of age $x-1$) at generation $n-1$, and
the same number $F_{n-2}^{x-2}$ (of age $x-2$) at generation $n-2$ and so on, up to generation $n-x+1$, in
which the rabbits were newborns. The condition  $x\le\min\{d,n\}$ ensures that $F_n^x$ is defined and that
$n-x+1\ge 1$ refers to an existing generation.
\end{proof}

\begin{lemma}
For $n\ge2$, the number of rabbits that die at generation $n$ is $F_{n-1}^d$.
\end{lemma}

\begin{proof}
Immediate from Definition~\ref{def:F_n^x}: $F_{n-1}^d$ is the number of
rabbits with age $d$. The condition $n\ge 2$ ensures that $F_{n-1}^d$ is defined.
\end{proof}

\begin{lemma}\label{lem:newborns}
The number $F_n^1$ of newborn rabbits at generation $n$ is equal to $\sum_{x=f}^{d} F_{n-1}^x$.
\end{lemma}

\begin{proof}
The number of newborns is equal to the number of fertile rabbits in the previous generation. Thus we need to consider
all the rabbits that in the previous generation have age at least $f$. By Definition~\ref{def:F_n^x}
we have that there are $F_{n-1}^f$ rabbits of age $f$, $F_{n-1}^{f+1}$ rabbits of age $f+1$, and
so on up to $F_{n-1}^d$ rabbits of age $d$. Thus the total number of rabbits that have a fertile
age, that is an age $\ge f$,  is $F_{n-1}^f+F_{n-1}^{f+1}+\ldots+F_{n-1}^d$.
\end{proof}

\bigskip

Now, armed with the above equalities, we can easily unravel the base equation. We distinguish four subcases.

\subsection{Case 1: $2\le n\le f$}

This case is trivial, but for the sake of clarity of exposition we treat it as the other cases.
In the first $f$ generations there are no newborns nor deaths. Hence we have
$newborns_n=deaths_n=0$ and since $F_1=1$, the base Equation~(\ref{eq:base})
becomes

\begin{equation}\label{eq:case1}
F_n = 1, \ \ \ \ \cfor 2 \le n \le d.
\end{equation}

\subsection{Case 2: $f+1\le n\le d$}

Since $n\le d$ there are no deaths, thus $deaths_n=0$. However some rabbits start to proliferate.
By Lemma~\ref{lem:newborns} we have that the newborns for generation $n$ are
$F_n^1=\sum_{x=f}^d F_{n-1}^x$.

By Lemma~\ref{lem:goback} we have that each element $F_{n-1}^x$ of the last summation can be
substituted by the equal term $F_{n-1-(f-1)}^{x-(f-1)}$, that we find by going up-left diagonally for $f-1$ generations,
and thus, starting from the base equation, we have that

\begin{eqnarray*}
F_n&=&F_{n-1}+newborns_n-deaths_n\\
      &=&F_{n-1}+F_{n-1}^1-0\\
      &=&F_{n-1}+\sum_{x=f}^{d-1} F_{n-1}^x\makebox{\ \ \ \ \ (by Lemma~\ref{lem:newborns})}\\
      &=&F_{n-1}+\sum_{x=f}^{d-1} F_{n-f}^{x-f+1}\makebox{\ \ \ \ \ (by Lemma~\ref{lem:goback})}\\
      &=&F_{n-1}+\sum_{x=1}^{d-f} F_{n-f}^{x}\makebox{\ \ \ \ \ (index substitution)}\\
      &=&F_{n-1}+\sum_{x=1}^{d} F_{n-f}^{x}\makebox{\ \ \ \ \ (added elements are 0)}\\
      &=&F_{n-1}+F_{n-f}\makebox{\ \ \ \ \ (by Equation~(\ref{eq:Fn=sum_of_row}))}\\
\end{eqnarray*}

Where the second to last  steps is true because for generation $n-f$ there are no rabbits with age $x>d-f$ since
the condition $n\le d$ implies that $n-f\le d-f$ and thus no rabbit can have an age bigger than
this number. Thus $F_{n-1}^x=0$ for $x>d-f$. 

Hence we have 
\begin{equation}\label{eq:case2}
F_n = F_{n-1}+F_{n-f}, \ \ \ \ \cfor f < n \le d.
\end{equation}

\subsection{Case 3: $n=d+1$}

This case is very similar to the previous one with the unique exception that we need to account for the first death:
indeed in generation $d$ the first rabbit dies, and it is the only one that dies, thus we have that $deaths_{n}=F_{n-1}^d=F_d^d=1$. The analysis of the newborns
carried out for the previous case applies also to this case. Hence we have that

\begin{equation}\label{eq:case3}
F_{d+1}=F_{d}-F_{d-f}-1.
\end{equation}

\subsection{Case 4:  $n\ge d+2$}

Next, we unravel Equation~(\ref{eq:base}) for the other values of $n$. The reasoning that we will provide
for this case cannot be applied to the previous ones because it involves $F_{n-d-1}$ which, for $n<d+2$,
is not defined.

\begin{eqnarray*}
F_n&=&F_{n-1}+newborns_n-deaths_n\\
      &=&F_{n-1}+F_n^1-F_{n-1}^d\\
      &=&F_{n-1}+\sum_{x=f}^{d} F_{n-1}^x-F_{n-1}^d\makebox{\ \ \ \ \ (by Lemma~\ref{lem:newborns})}\\
      &=&F_{n-1}+\sum_{x=f}^{d-1} F_{n-1}^x
\end{eqnarray*}

As done for case 2, by Lemma~\ref{lem:goback} we have that each element $F_{n-1}^x$ of the last summation can be
substituted by the equal term $F_{n-1-(f-1)}^{x-(f-1)}$, that we find by going up-left diagonally for $f-1$ generations,
and thus we have that

\begin{eqnarray*}
F_n&=&F_{n-1}+\sum_{x=f}^{d-1} F_{n-1}^x\\
      &=&F_{n-1}+\sum_{x=f}^{d-1} F_{n-f}^{x-f+1}\makebox{\ \ \ \ \ (by Lemma~\ref{lem:goback})}\\
      &=&F_{n-1}+\sum_{x=1}^{d-f} F_{n-f}^{x}\makebox{\ \ \ \ \ (index substitution)}\\
\end{eqnarray*}

We now observe that 
$$F_{n-f} = \sum_{x=1}^{d} F_{n-f}^x = \sum_{x=1}^{d-f} F_{n-f}^x+\sum_{x=d-f+1}^{d} F_{n-f}^x$$
and thus we have that
$$\sum_{x=1}^{d-f} F_{n-f}^x = F_{n-f} - \sum_{x=d-f+1}^{d} F_{n-f}^x.$$

\noindent
Hence, we have that
\begin{eqnarray*}
F_n & = & F_{n-1}+\sum_{x=1}^{d-f} F_{n-f}^{x}\\
       & = & F_{n-1} + F_{n-f} - \sum_{x=d-f+1}^{d} F_{n-f}^x.
\end{eqnarray*}

As done before we can now shift up-left for $d-f$ generations the terms of the summation using Lemma~\ref{lem:goback}:
\begin{eqnarray*}
\sum_{x=d-f+1}^{d} F_{n-f}^x & = & \sum_{x=d-f+1-(d-f)}^{d-(d-f)} F_{n-d}^x\makebox{\ \ \ \ \ (by Lemma~\ref{lem:goback})}\\
					     & = & \sum_{x=1}^{f} F_{n-d}^x\makebox{\ \ \ \ \ (index substitution)}\\
\end{eqnarray*}
Thus we have that
\begin{eqnarray}
\nonumber
F_n  & = & F_{n-1} + F_{n-f} - \sum_{x=d-f+1}^{d} F_{n-f}^x\\
\label{eq:one}
        & = & F_{n-1} + F_{n-f} -\sum_{x=1}^{f} F_{n-d}^x.
\end{eqnarray}

%%%%
And to conclude the unraveling of the formula we observe that the summation of the negative terms $F_{n-d}^x$
is equal to the total number of rabbits $F_{n-d-1}$ in the previous generation.
Indeed the first term $F_{n-d}^1$ of this summation is the number of newborns of generation $n-d$ which are equal
to the number of fertile rabbits in generation $n-d-1$; that is
$$F_{n-d}^1=\sum_{x=f}^d F_{n-d-1}^x$$
and the other terms, shifting them up-left of 1 generation, using again Lemma~\ref{lem:goback}, are equal to
$$\sum_{x=2}^{f} F_{n-d}^x = \sum_{x=1}^{f-1} F_{n-d-1}^x.$$

By putting together these facts we have that
\begin{equation}\label{eq:two}
\sum_{x=1}^{f} F_{n-d}^x = F_{n-d}^1+\sum_{x=2}^{f} F_{n-d}^x =
\sum_{x=f}^d F_{n-d-1}^x + \sum_{x=1}^{f-1} F_{n-d-1}^x = F_{n-d-1}.
\end{equation}

\bigskip
\noindent

Putting together Equations~(\ref{eq:one}) and (\ref{eq:two}) we have that
\begin{equation}\label{eq:case4}
F_n = F_{n-1}+F_{n-f}-F_{n-d-1}, \ \ \ \ \cfor n\ge d+2.
\end{equation}

\subsection{The formula}

To summarize, we provide the following theorem.

\begin{thm}\label{thm:formula}
Let $F_1=1$ and let $f$ and $d$ be integers such that $1\le f \le d$. The number $F_n$ of rabbits at generation $n$, 
for a population of rabbits that become fertile at age $f$ and die
at age $d\ge f$, is given by
\begin{equation}
F_n = 
\begin{cases}
1,                                                       & \cfor 2\le n\le f\\
F_{n-1}+F_{n-f},                                & \cfor f < n \le d\\
F_{n-1}+F_{n-f}-1,                             & \cfor n = d+1\\
F_{n-1}+F_{n-f}-F_{n-d-1},                & \cfor n \ge d+2\\
\end{cases}
\end{equation}
\end{thm}

\begin{proof}
From Equations~(\ref{eq:case1}), (\ref{eq:case2}), (\ref{eq:case3}), and (\ref{eq:case4}).
\end{proof}

The Fibonacci sequence is obtained  with $f=2$ and $d=\infty$.
Using $f=2$ and $d=3$ one obtains the Padovan sequence.

\section{Conclusions}

We have given a simple formula, proved with a straightforward argumentation, for the number of rabbits in the generalized Fibonacci problem, in which rabbits
become fertile after an arbitrary number of generations and they also, at some point, die. The dying rabbits problem
was posed by Brother U. Alfred in the first issue of the Fibonacci Quarterly~\cite{A1963-1}, probably as what the 
author expected to be an easy counting problem. Here is a verbatim quote from a subsequent paper~\cite{A1963-4} by Brother U. Alfred:
{\em Originally, it was thought that the rabbits removed would constitute a sequence which could be readily identified 
with an expression involving Fibonacci numbers. But after several attempts by a number of people it appeared
that it would be difficult to arrive at an answer by straightforward intuition.}
In this paper we have shown that, after all, Brother U. Alfred was right in considering the dying rabbits problem
an easy counting problem.

\bibliographystyle{plain}

\end{document}